\documentclass[11pt]{article}
\usepackage{amsmath, amsthm,  latexsym, amssymb, amsfonts, epsfig, graphicx}
\usepackage{mathrsfs}
\allowdisplaybreaks

\renewcommand{\d}{\mathrm{d}}
\newcommand{\esssup}{\mathop{\mathrm{esssup}}}
\newcommand{\ep}{\varepsilon}

\newcommand{\R}{\mathbb{R}}

\newcommand{\E}{\mathbb{E}}
\newcommand{\PP}{\mathbb{P}}
\newcommand{\T}{\mathcal{T}}
\newcommand{\K}{\mathcal{K}}

\newcommand{\F}{\mathcal{F}}
\newcommand{\lt}{\left}
\newcommand{\rt}{\right}
\newtheorem{theorem}{Theorem}[section]

\newtheorem{lemma}[theorem]{Lemma}
\par
\newtheorem{definition}[theorem]
{Definition}
\par

\par

\par
\newtheorem{remark}[theorem]{Remark}
\par

\par

\par

\begin{document}

\title{Peng's Maximum Principle for a Stochastic Control Problem Driven by a Fractional and a Standard Brownian Motion}
\author{BUCKDAHN Rainer$^{a,b}$\ and
JING Shuai$^c$\footnote{Corresponding author. E-mail: shuaijingsj@gmail.com. Supported by National Natural Science Foundation of China (Project No. 11301560).}\\
\small{$^a$ D\'epartement de Math\'ematiques, Universit\'e de Bretagne Occidentale, 29285, Brest, France}\\
\small{$^b$ School of Mathematics, Shandong University, 250100, Jinan, P.R.China}\\
\small{$^c$ School of Management Science and Engineering, Central University
of Finance and Economics, }\\
\small{100081,
Beijing, P.R.China}}

\date{}\maketitle

\begin{abstract}
We study a stochastic control system involving both a standard and a fractional Brownian motion with Hurst parameter less than 1/2. We apply an anticipative Girsanov transformation to transform the system into another one, driven only by the standard Brownian motion with coefficients depending on both the fractional Brownian motion and the standard Brownian motion. We derive a maximum principle and the associated stochastic variational inequality, which both are generalizations of the classical case.
\end{abstract}
\medskip

\medskip


\section{Introduction}
\label{sec:1}
We study a control problem which controlled state process is driven by both a standard Brownian motion and a fractional Brownian motion with Hurst parameter $H\in(0,1/2)$, and we derive the stochastic maximum principle and the associated variational inequality. To be more precise, we consider the state process governed by the following  controlled stochastic differential equation
\begin{equation}\label{eq:state}
\lt\{
\begin{array}{l}
\d X^u(t)=\sigma(t) X^u(t)\d B^H(t) + \beta(t,X^u(t),u(t))\d W(t)+ b(t, X^u(t), u(t))\d t, \quad t\in[0,T],\\
X^u(0)=x_0,
\end{array}\rt.
\end{equation}
where the functions $\sigma$, $\beta$ and $b$ are introduced in Section 2, and the control process $u$ takes values in a metric space $U$. Thus, in our framework, the diffusion part consists of two parts: one is represented by a stochastic integral with respect to the fractional Brownian motion $B^H$, which integrand is linear in the state process, and the other by an It\^o integral with respect to the Brownian motion $W$, which integrand is nonlinear in the state process. The control problem consists in minimizing the cost functional defined as follows
\begin{equation*}
J(u)=\E\lt[\Phi(X^u(T))+\int^T_0 f(t,X^u(t),u(t))\d t\rt],
\end{equation*}
where the functions $\Phi$ and $f$ are introduced in Section 2.

Stochastic differential equations driven by both a standard Brownian motion and a fractional Brownian motion have been studied by several authors, for example, for the case $H\in(1/2, 1)$ by Guerra and Nualart \cite{GuNu}, Mishura and Shevchenko \cite{MiSh}. The properties of the fractional Brownian motion with Hurst parameter $H\in(1/2, 1)$ and $H\in(0, 1/2)$ are quite different. Besides a pathwise definition of the integral, the classical divergence operator is widely used to define the stochastic integral when $H\in(1/2, 1)$. However, in the case $H\in(0,1/2)$, the domain of the divergence operator becomes too small. For instance, Cheridito and Nualart \cite{ChNu} showed that even the fractional Brownian motion itself is not included in the domain. To overcome this difficulty, Cheridito and Nualart \cite{ChNu} and Le\'on and Nualart \cite{LeNu} defined a new type of operator and called it the extended divergence operator. By using the extended divergence operator, Le\'on and San Mart\'in \cite{LeSa} studied linear stochastic differential equations driven by a fractional Brownian motion with $H\in(0,1/2)$ with the help of the chaos decomposition. Jien and Ma \cite{JiMa} worked on stochastic differential equations driven by fractional Brownian motions by applying anticipative Girsanov transformations developed by Buckdahn \cite{Bu}, while Jing and Le\'on \cite{JL} also made use of the Girsanov transformation method to deal with semilinear backward doubly stochastic differential equations driven by a Brownian motion and a fractional Brownian motion with $H\in(0,1/2)$ and the associated stochastic partial differential equations driven by the fractional Brownian motion.

The variational inequality and stochastic maximum principle for controlled systems driven by a Brownian motion have been investigated by many authors. Without being exhaustive, let us mention among them, for instance, Bismut \cite{Bi}, Bensoussan \cite{Be}, Peng \cite{P90} and Buckdahn \textit{et al.} \cite{BuDjLi}. However, for controlled systems involving fractional Brownian motions, there are only very few works. Biagini \textit{et al.} \cite{BHOS} studied a stochastic maximum principle for processes driven only by an $m$-dimensional fractional Brownian motion with $H\in(1/2,1)^m$ and derived an adjoint linear fractional backward stochastic differential equation.  Hu and Zhou \cite{HuZh} considered an optimal control problem of stochastic linear systems involving a fractional Brownian motion with Hurst parameter $H\in(0,1/2)$, and they introduced a Riccati equation which is a backward stochastic differential equation driven by the fractional Brownian motion and a classical Brownian motion. It is worth noting that this Brownian motion is the one that generates the fractional Brownian motion, hence they are not independent. Han \textit{et al.} \cite{HaHuSo} obtained a stochastic maximum principle for a stochastic control problem defined through a general controlled system driven by a fractional Brownian motion with $H>1/2$. Similar to \cite{HuZh}, their adjoint backward stochastic differential equation is driven by the fractional Brownian motion and its underlying Brownian motion.

Here, in our framework, the controlled system involves both a standard and a fractional Brownian motion with $H\in(0,1/2)$. We use the extended divergence operator to define the stochastic integral with respect to the fractional Brownian motion. The linearity of the integrand of the integral with respect to $B^H$ in the state process allows, similarly to \cite{JL}, to apply the anticipative Girsanov transformation  to transform the original controlled system into another one driven only by the standard Brownian motion $W$, but with coefficients depending on the paths of both $W$ and $B^H$. Our adjoint backward stochastic differential equation involves, besides the Brownian martingale, also an orthogonal martingale, which is a Brownian martingale in the classical case (see Section 4). This orthogonal martingale comes from the Galtchouk-Kunita-Watanabe decomposition. Such backward stochastic differential equations were employed by Buckdahn and Ichihara \cite{BuIc} and Buckdahn \textit{et al.} \cite{BuLaRaTa} to study optimal control systems and associated Hamilton-Jacobi-Bellman equations. In our work here, we compare our main result with the classical characterization of an optimal control and we show that, if we replace the fractional Brownian motion with a standard Brownian motion, i.e., if we apply our Girsanov transformation in the classical, Brownian framework, we get the same result. Hence, our result indeed generalizes the classical one.

In this paper we deal only with the case $H\in (0,1/2)$ since we use the extended divergence operator as stochastic integral with respect to the fractional Brownian motion. Nevertheless, when using the divergence operator in the case $H\in(1/2,1)$, our method is still valid and the computations are even easier. The key difference between the two cases relies mainly on the distinct definitions of the divergence operator and the extended divergence operator.

The paper is organized as follows: In Section 2 we recall some preliminaries, i.e., some basic settings and some basics on the fractional Brownian motion, the extended divergence operator and the Girsanov transformation. Our main results,  the variational inequality and the stochastic maximum principle, are stated in Section 3. Finally, in Section 4 we compare our result with Peng's criterion for the optimality of a stochastic control in the Brownian setting \cite{P90}. The proofs of the results in Section 3 are given in the Appendix to improve the readability.
\section{Preliminaries}
\subsection{General Setting and Fractional Brownian Motion}
Let $T>0$ be a fixed time horizon. Let $\{W(s), s\in[0,T]\}$ be a standard Brownian motion
on a complete probability space $(\Omega_1, \F_1, \PP_1)$ and $\{B^H(s), s\in[0,T]\}$ be a
fractional Brownian motion with Hurst parameter $H\in(0,1/2)$ defined
on another complete probability space  $(\Omega_2, \F_2,$ $\PP_2)$.
We introduce $(\Omega, \F, \PP)$ as the product space $(\Omega,\F,\PP) =
(\Omega_1,\F_1,\PP_1) \otimes(\Omega_2,\F_2,\PP_2)
=(\Omega_1\times\Omega_2,\F_1\otimes\F_2,\PP_1\otimes\PP_2)$ which we suppose to be completed.
The processes $W$ and $B^H$ are canonically extended from
$(\Omega_1,\F_1,$ $\PP_1)$ and $(\Omega_2, \F_2,$ $\PP_2)$, respectively,
to the product space $(\Omega, \F, \PP)$.

We define three filtrations: one is generated by the Brownian motion: $\mathbb{F}^W=\{\mathcal{F}_t^W = \sigma\{W(s), 0\leq s\leq
t\}\vee{\cal N}, t\in[0,T]\}$, one is generated by the fractional Brownian motion: $\mathbb{F}^B=\{\mathcal{F}_t^{B}=\sigma\{B^H(s), 0\leq s\leq
t\}\vee{\cal N}, t\in[0,T]\},$  and another one is generated by the Brownian motion $W$ and  the fractional Brownian motion $B^H$ over the time interval $[0,T]$ : $\mathbb{H}=\{\mathcal {H}_t
=\mathcal{F}_{t}^W \vee \mathcal{F}_t^{B}\vee \mathcal{N}, {t\in[0,T]}\}$. Here $\mathcal{N}$ denotes the set of all $\PP$-null sets.

For $p>1$, we denote by $L^p_{\mathbb{H}}(\Omega\times[0,T])$ the space of real valued $\mathbb{H}$-adapted processes such that
$$
\|\varphi\|_{L^p_{\mathbb{H}}}=\lt(\E\lt[\int^T_0\lt|\varphi(t)\rt|^p\d t\rt]\rt)^{1/p}<+\infty.
$$

It is well-known that, for $H\in(0,1/2)$, there exists another canonical Wiener process $W^0$ on $(\Omega_2, \F_2,$ $\PP_2)$ such that we have the following representation:
\[
B^H(t)=\int^t_0 K_H(t,s) \d W^0(s),
\]
where
\[
K_H(t,s)= C_H \lt[\lt(\frac{t}{s}\rt)^{H-1/2}t^{H-1/2}-(H-1/2)s^{1/2-H}\int^t_s u^{H-3/2}(u-s)^{H-1/2}\d u\rt],
\]
and
\[
C_H=\sqrt{\frac{2H}{(1-2H)\beta(1-2H,H+1/2)}}.
\]
Hence, the process $B^H$ is a centered Gaussian process with covariance function
\[
R_H(t,s)=\E\lt[B^H(t)B^H(s)\rt]=\frac12\lt(t^{2H}+s^{2H}-|t-s|^{2H}\rt).
\]
\subsection{Extended Divergence Operator}
We briefly recall the definition of the extended divergence operator as the stochastic integral with respect to the fractional Brownian motion $B^H$; for more details, we refer to \cite{JL}. The extended divergence operator was first studied by Cheridito and Nualart \cite{ChNu} and further investigated by Le\'on and Nualart \cite{LeNu}.

To this end, we define a Hilbert  space $\mathcal{H}_H$ as the completion of the space of step functions over $[0,T]$ with respect to the inner product
\[
\langle I_{[0,t]}, I_{[0,s]}\rangle_{\mathcal{H}_H}=R_H(t,s), \ t,s\in[0,T].
\]
On the space $\mathcal{H}_H$, an isometry $\mathcal{H}_H\ni\varphi\to B^H(\varphi)\in L^2(\Omega,\F,\PP)$ is defined by extending the map $I_{[0,t]}\to B^H_t$. Moreover, by the transfer principle (see Nualart \cite{Nu}), one has the existence of an operator $\K:\mathcal{H}_H\to L^2([0,T])$ such that
\[
B^H(\varphi)=\int^T_0 (\K \varphi)(s)\d W^0 (s), \ \varphi\in\mathcal{H}_H,\ \ \textrm{and}\ \ (\K I_{[0,t]})(s)=K_H(t,s), \ s,t\in[0,T].
\]
We denote by $\K^\ast$ its adjoint operator.

Let $\mathcal{S}_{\K}$ be the class of all smooth functionals of the form
\[
F=f\lt(B^H(\varphi_1),\cdots, B^H(\varphi_m), W(\psi_1),\cdots, W(\psi_n)\rt),\ \ m,n\ge 1,
\]
where $\varphi_1,\cdots,\varphi_m$ are elements of $\mathcal{H}_H$, $\psi_1,\cdots,\psi_n\in L^2([0,T])$, $W(\psi_1),\cdots, W(\psi_n)$ are Wiener integrals of $\psi_1,\cdots,\psi_n$ with respect to $W$, and $f\in C^{\infty}_p(\R^{m+n})$  - the space of all $C^{\infty}$ function over $\R^{m+n}$, which together with all their derivatives are of polynomial growth.

A smooth functional $F\in\mathcal{S}_{\K}$ of above form has Malliavin derivatives with respect to $B^H$ and $W$ defined as follows:
\[
D^BF=\sum^m_{i=1}\frac{\partial f}{\partial x_i}\lt(B^H(\varphi_1),\cdots, B^H(\varphi_m), W(\psi_1),\cdots, W(\psi_n)\rt)\varphi_i,
\]
and
\[
D^WF=\sum^n_{i=1}\frac{\partial f}{\partial x_{m+i}}\lt(B^H(\varphi_1),\cdots, B^H(\varphi_m), W(\psi_1),\cdots, W(\psi_n)\rt)\psi_i.
\]
We remark that both $D^B F$ and $D^W F$ are in $L^p(\Omega\times[0,T]; \mathcal{H}_H)$, for all $p\ge 2$.

For $u\in L^2(\Omega\times[0,T])$, we define the following stochastic integrals with respect to $B^H$ and $W$, respectively.

\begin{definition}
Let $u\in L^2(\Omega\times[0,T])$.
If there exists a random variable $\delta^B(u)\in L^2(\Omega,\F,\PP)$ such that
\begin{equation}\label{eq_defdeltab}
\E\lt[\langle \mathcal{K}^\ast\mathcal{K}D^BF,u\rangle_{L^2{([0,T])}}\rt]=\E\lt[F\delta^B(u)\rt], \ \ \textrm{for all}\ \ F\in\mathcal{S}_\K,
\end{equation}
we say $u\in Dom\ \delta^B$ and call $\delta^B(u)$ the extended divergence operator of $u$ with respect to $B^H$.
\end{definition}

\begin{definition}
Let $u\in L^2(\Omega\times[0,T])$.
If there exists a random variable $\delta^W(u)\in L^2(\Omega,\F,\PP)$ such that
\begin{equation}\label{eq_defdeltaw}
\E\lt[\langle D^WF,u\rangle_{L^2{([0,T])}}\rt]=\E\lt[F\delta^W(u)\rt], \ \ \textrm{for all}\ \ F\in\mathcal{S}_\K,
\end{equation}
we say $u\in Dom\ \delta^W$ and call $\delta^W(u)$ the Skorohod integral of $u$ with respect to $W$.
\end{definition}

\begin{remark}\label{remark}
1. Given $u\in L^2(\Omega\times[0,T])$ and $t\in[0,T]$ such that $uI_{[0,t]}\in Dom \ \delta^B$, we write $\int^t_0u(s)\d B^H(s)$ for $\delta^B(uI_{[0,t]})$.

2. If $u\in L^2(\Omega\times[0,T])$ is $\mathbb{H}$-adapted, then the Skorohod integral $\delta^W(u)$ exists and it coincides with the It\^o integral $\int^T_0u(s)\d W(s)$ (Recall that $W$ is an $\mathbb{H}$-Brownian motion).
\end{remark}

\subsection{Girsanov Transformations}
The Girsanov transformation with respect to the fractional Brownian motion constitutes an essential tool in our approach for our stochastic control problem.

Throughout this paper we use the following hypothesis.
\medskip

\textbf{(H1)} Let $\sigma:[0,T]\to\mathbb{R}$ be a square integrable Borel function such that $\sigma I_{[0,t]}$ belongs to $\mathcal{H}_H$, for every $t\in[0,T]$, and $\sup_{0\le t\le T}\int^t_0((\mathcal{K}\sigma I_{[0,t]})(r))^2\d r<+\infty.$
\medskip

Recall that hypothesis \textbf{(H1)} is in particular satisfied if $\sigma(t)=\sigma, t\in[0,T]$, for some constant $\sigma\in\R$, and $(\mathcal{K}\sigma I_{[0,t]})(r)=K_H(t,r)$, $(t,r)\in[0,T]^2$.

For $t\in[0,T]$, we consider the following transformations on $\Omega_2$:
$$
\mathcal{T}_t(\omega_2)=\omega_2+\int^{\cdot\wedge t}_0(\mathcal{K}\sigma I_{[0,t]})(r)\d r,
\ \ \omega_2\in\Omega_2,\ \ t\in[0,T],
$$
and
$$
\mathcal{A}_t(\omega_2)=\omega_2-\int^{\cdot\wedge t}_0(\mathcal{K}\sigma I_{[0,t]})(r)\d r,
\ \ \omega_2\in\Omega_2,\ \ t\in[0,T].
$$
The Girsanov Theorem (see for example Buckdahn \cite{Bu}) gives
that for any square integrable random variable $F$, we have
\begin{equation}\label{eq:gir}
\E[F]=\E[F(\mathcal{A}_t)\kappa_t]=\E[F(\mathcal{T}_t)\kappa_t^{-1}(\mathcal{T}_t)],
\end{equation}
where
\begin{equation}\label{ep}
\kappa_t=\exp\lt(\int^t_0\sigma(r)\d B^H(r)-
\frac12\int^t_0((\mathcal{K}\sigma I_{[0,t]})(r))^2\d r\rt).
\end{equation}
From Lemma 2.4 in \cite{JL} we have that
\begin{equation}\label{eq:ep}
E\lt[\sup_{0\le t\le T}\kappa^p_t\rt]<+\infty, \ \
E\lt[\sup_{0\le t\le T}\kappa^p_t(\mathcal{T}_t)\rt]<+\infty,
\ \mbox{for\ all}\ p\in\R.
\end{equation}

\section{Variational Inequality and the Maximum Principle}

\subsection{The Stochastic Control Problem}

Let $U$ be a nonempty subset of $\R^k$. Let
$\{u(s), 0\le s\le T\}$ be an admissible control process,
which takes values in $U$ and is  $\mathbb{H}$-adapted,
such that
$$
\esssup_{0\le t\le T} \E[|u(t)|^p]<+\infty, \quad \textrm{for}\ \ p\ge1.
$$
The set of admissible control processes is denoted by
$\mathcal{U}_{ad}$. From \textbf{(H1)} and (\ref{eq:ep})
we get that if $\{u(s), 0\le s\le T\}$ is an admissible control, then both
$\{u(s,\mathcal{T}_s), 0\le s\le T\}$ and $\{u(s,\mathcal{A}_s), 0\le s\le T\}$
are admissible controls. In particular, we have
\begin{equation}\label{eq:control}
\esssup_{0\le t\le T} \E[|u(t,\mathcal{T}_t)|^p]<+\infty \ \ \textrm{and}\ \
\esssup_{0\le t\le T} \E[|u(t,\mathcal{A}_t)|^p]<+\infty, \quad p\ge1.
\end{equation}

We consider the following stochastic control system:
\begin{equation}\label{eq_sys}
\lt\{
\begin{array}{l}
\d X^u(t)=\sigma(t) X^u(t)\d B^H(t) + \beta(t,X^u(t),u(t))\d W(t)+ b(t, X^u(t), u(t))\d t, \quad t\in[0,T],\\
X^u(0)=x_0.
\end{array}\rt.
\end{equation}
Notice that only the coefficients $\beta$ and $b$ depend on the control, but not $\sigma$.
Moreover, the stochastic integral with respect to the fractional Brownian motion is linear in $X^u$
and is interpreted in the extended divergence sense.
The cost functional is defined by
\begin{equation}\label{eq_costf}
J(u)=\E\lt[\Phi(X^u(T))+\int^T_0 f(t,X^u(t),u(t))\d t\rt].
\end{equation}
Our control problem consists in minimizing the cost functional  $J(u)$ over $\mathcal{U}_{ad}$.

Now we state the assumptions on the coefficients:
$$
\beta, b, f: [0,T]\times\R\times\R^k\to\R,\ \ \Phi: \R\to\R.
$$

\textbf{(H2)} The functions $\beta, b, f, \Phi$ are twice differentiable with respect to $x$.
Moreover, $\beta, b, f, \Phi$ and their derivatives
$\beta_x, b_x, f_x, \Phi_x$ $\beta_{xx}, b_{xx}, f_{xx}, \Phi_{xx}$
are continuous in $(x,u)$ and bounded, uniformly with respect to $(t,u)\in[0,T]\times U$.

\subsection{Main Results}
In this subsection we state our main results, i.e., the variational inequality and the maximum principle.

First we state the following important theorem. It helps us to establish a link between the semilinear stochastic differential equation (\ref{eq_sys}), driven by both the standard Brownian motion $W$ and the fractional Brownian motion $B^H$, and a stochastic differential equation driven only by the standard Brownian motion $W$, with coefficients depending on the fractional Brownian motion.
\begin{theorem}\label{thm_girsanov}
The process $X^u:=\lt\{ X^u(t)=\zeta^u(t,\mathcal{A}_t)\kappa_t, t\in[0,T]\rt\}$ is the unique solution of equation (\ref{eq_sys}) in $L^2_{\mathbb{H}}(\Omega\times[0,T])$,  where $\zeta^u$ is the unique solution of the pathwise stochastic differential equation
\begin{equation}\label{eq_zeta}
\left\{
\begin{array}{l}
\d \zeta^u(t)=\kappa_t^{-1}(\mathcal{T}_t) \beta(t,\zeta^u(t)\kappa_t(\mathcal{T}_t),u(t,\mathcal{T}_t))\d W(t)+ \kappa_t^{-1}(\mathcal{T}_t)b(t,\zeta^u(t)\kappa_t(\mathcal{T}_t), u(t,\mathcal{T}_t))\d t, \ t\in[0,T];\\
\zeta^u(0)=x_0.
\end{array}
\right.
\end{equation}
\end{theorem}

For the reader's convenience we give the proof; it is shifted to the Appendix.
\medskip

The above theorem allows to rewrite the cost functional (\ref{eq_costf}) as follows:
\begin{equation}\label{eq_cost}
J(u)=\E\lt[\Phi\big(\zeta^u(T)\kappa_T(\mathcal{T}_T)\big)\kappa_T^{-1}(\mathcal{T}_T)
+\int^T_0f\big(t,\zeta^u(t)\kappa_t(\mathcal{T}_t),u(t,\mathcal{T}_t)\big)\kappa_t^{-1}(\mathcal{T}_t)\d t\rt].
\end{equation}
We have transformed our stochastic control problem into a formally classical control problem which contains the fractional Brownian motion implicitly.

Since the control process $u(t,\mathcal{T}_t)$ appearing in (\ref{eq_zeta}) and (\ref{eq_cost}) contains always the transformation $\mathcal{T}_t$, for the simplicity of notations, we denote it by $v(t)$, i.e., $v(t)=u(t,\mathcal{T}_t)$. From (\ref{eq:control}) we know that both $u$ and $v$ are admissible controls.

Let us now suppose that $(y(\cdot), v(\cdot))$ is an optimal solution of the control problem, i.e.,
\begin{equation}\label{eq_y}
\left\{
\begin{array}{l}
\d y(t)=\kappa_t^{-1}(\mathcal{T}_t) \beta(t,y(t)\kappa_t(\mathcal{T}_t),v(t))\d W(t)+ \kappa_t^{-1}(\mathcal{T}_t)b(t,y(t)\kappa_t(\mathcal{T}_t), v(t))\d t, \ t\in[0,T],\\
y(0)=x_0,
\end{array}
\right.
\end{equation}
and
$$
J(v)=\inf_{u\in\mathcal{U}_{ad}} J(u).
$$
Following Peng's approach \cite{P90}, we construct a perturbed admissible control as follows:
$$
v^\ep(t)=
\lt\{
\begin{array}{ll}
\tilde{v}(t) & \tau-\ep\le t\le \tau+\ep,\\
v(t) & \textrm{otherwise},
\end{array}
\rt.
$$
where $0<\tau<T$ is arbitrarily fixed, $\ep>0$ is arbitrarily chosen such that $[\tau-\ep,\tau+\ep]\subset [0, T]$, and $\tilde{v}$ is an arbitrary bounded admissible control from $\mathcal{U}_{ad}$. Let $y^{\ep}(\cdot)$ be the solution of (\ref{eq_y}) with $v^\ep$ at the place of $v$. Then from the setting of the control problem, we have
$$J(v^\ep)-J(u)\ge 0.$$

Let $y_1(\cdot)$ and $y_2(\cdot)\in L^{\infty,-}_{\mathbb{H}}([0,T]) (:=\bigcap_{p\ge2}L^p_{\mathbb{H}}(\Omega\times[0,T]))$ be the solutions of the equations
\begin{equation}\label{eq_y_1}
\begin{aligned}
y_1(t)=&\int^t_0\bigg[b_x(s,y(s)\kappa_s(\mathcal{T}_s),v(s))y_1(s)\\
&\qquad+\kappa_s^{-1}(\mathcal{T}_s)
\bigg(b(s,y(s)\kappa_s(\mathcal{T}_s),v^\ep(s))-b(s,y(s)\kappa_s(\mathcal{T}_s),v(s))\bigg)\bigg]\d s\\
&+\int^t_0\bigg[\beta_x(s,y(s)\kappa_s(\mathcal{T}_s),v(s))y_1(s)\\
&\qquad\quad+\kappa_s^{-1}(\mathcal{T}_s)
\bigg(\beta(s,y(s)\kappa_s(\mathcal{T}_s),v^\ep(s))-\beta(s,y(s)\kappa_s(\mathcal{T}_s),v(s))\bigg)\bigg]\d W(s),
\end{aligned}
\end{equation}
and\begin{equation}\label{eq_y_2}
\begin{aligned}
y_2(t)=&\int^t_0\bigg[b_x(s,y(s)\kappa_s(\mathcal{T}_s),v(s))y_2(s)+\frac12\kappa_s(\mathcal{T}_s)
b_{xx}(s,y(s)\kappa_s(\mathcal{T}_s),v^\ep(s))y_1^2(s)\bigg]\d s\\
&+\int^t_0\bigg[\beta_x(s,y(s)\kappa_s(\mathcal{T}_s),v(s))y_2(s)+\frac{1}{2}\kappa_s(\mathcal{T}_s)
\beta_{xx}(s,y(s)\kappa_s(\mathcal{T}_s),v^\ep(s))y_1^2(s)\bigg]\d W(s)\\
&+\int^t_0\bigg(b_x(s,y(s)\kappa_s(\mathcal{T}_s),v^\ep(s))-b_x(s,y(s)\kappa_s(\mathcal{T}_s),v(s))\bigg)y_1(s)\d s\\
&+\int^t_0\bigg(\beta_x(s,y(s)\kappa_s(\mathcal{T}_s),v^\ep(s))-\beta_x(s,y(s)\kappa_s(\mathcal{T}_s),v(s))\bigg)y_1(s)\d W(s).
\end{aligned}
\end{equation}

We need the following estimates for $y_1$ and $y_2$.

\begin{lemma}\label{lemma_y1y2}
Under our hypotheses \textbf{(H1)} and \textbf{(H2)}, for any $p\ge2$, there is some $C_p\in\R_+$ independent of $\ep$ such that
\begin{equation}\label{est_y1}
\E\lt[\sup_{0\le t\le T}|y_1(t)|^{p}\rt]\le C_p\ep^{p/2},
\end{equation}
\begin{equation}\label{est_y2}
\E\lt[\sup_{0\le t\le T}|y_2(t)|^p\rt]\le C_p\ep^p.
\end{equation}
\end{lemma}
\noindent\textit {Proof:} First we prove inequality (\ref{est_y1}). From equation (\ref{eq_y_1}), using \textbf{(H2)} and the Buckh\"older-Davis-Gundy inequality, we obtain that, for $p\ge 2$,
\begin{equation}
\begin{aligned}
\E\lt[\sup_{0\le r\le t}|y_1(r)|^p\rt]&
\le C_p\E\lt[\sup_{0\le r\le t}\int^r_0|y_1(s)|^p\d s\rt]+ C_p\E\lt[\lt(\int^t_0I_{[\tau-\ep, \tau+\ep]}\kappa_s^{-2}(\mathcal{T}_s)\d s\rt)^{p/2}\rt]\\
&\le C_p\E\lt[\int^t_0\sup_{0\le r\le s}|y_1(r)|^p\d s\rt]+ C_p\ep^{p/2}\E\lt[\sup_{0\le s\le T}\kappa_s^{-p}(\mathcal{T}_s)\rt]\\
&\le C_p\E\lt[\int^t_0\sup_{0\le r\le s}|y_1(r)|^p\d s\rt]+ C_p\ep^{p/2},
\end{aligned}
\end{equation}
where the constant $C_p$ can be chosen independent of $t$. By the Gronwall inequality we get that
$$
\E\lt[\sup_{0\le t\le T} |y_1(t)|^p\rt]\le C_p\ep^{p/2}.
$$

Now we prove (\ref{est_y2}). From equation (\ref{eq_y_2}) and \textbf{(H2)}, applying Cauchy-Schwarz inequality and Buckh\"older-Davis-Gundy inequality, we have
\[
\begin{aligned}
\E\lt[\sup_{0\le r\le t}|y_2(r)|^p\rt]
\le& C_p \E\lt[\int^t_0\lt(\sup_{0\le r\le s}|y_2(r)|^p+\sup_{0\le r\le s}\kappa_r^p(\T_r)|y_1(r)|^{2p}\rt)\d s\rt]\\
&+C_p\E\lt[\lt|\int^t_0|I_{[\tau-\ep,\tau+\ep]}(s)y_1(s)|\d s\rt|^p\rt]+C_p\E\lt[\int^t_0I_{[\tau-\ep,\tau+\ep]}(s)\lt|y_1(s)\rt|^2\d s\rt]^{\frac{p}{2}}\\
\le&C_p \E\lt[\int^t_0\sup_{0\le r\le s}|y_2(r)|^p\d s\rt]+C_p\lt(\E\lt[\sup_{0\le r\le T}\kappa_r^{2p}(\T_r)\rt]\rt)^{\frac12}\lt(\E\lt[\sup_{0\le r\le T}y_1^{4p}(r)\rt]\rt)^{1/2}\\
&+C_p\ep^{p}\E\lt[\sup_{0\le t\le T}|y_1(t)|^p\rt]+C_p\ep^{p/2}\E\lt[\sup_{0\le t\le T}|y_1(t)|^p\rt].
\end{aligned}
\]
Hence, from (\ref{est_y1}) and the Gronwall inequality, we obtain
\begin{equation*}
\E\lt[\sup_{0\le t\le T}|y_2(t)|^p\rt]\le C_p\ep^p.
\end{equation*}
The proof is complete. $\Box$
\medskip

Set $y_3=y_1+y_2$. To derive our variational inequality, it is necessary to prove the following estimate.
\begin{lemma}\label{lemma_y-y3}
Under the hypothesis \textbf{(H2)}, for any $p\ge 2$, we have
\begin{equation}\label{eq_y1-y3}
\sup_{0\le t\le T}\E\lt[|y^{\ep}(t)-y(t)-y_3(t)|^p\rt]= o(\ep^p).
\end{equation}
\end{lemma}
For convenience of the reader the proof is given in the Appendix.
\medskip

The next lemma plays an important role in deriving the variational inequality.
\begin{lemma}\label{lemma_var}
Under the hypothesis \textbf{(H2)} we have
\begin{equation}\label{eq_var}
\begin{aligned}
&-\E\lt[\int^T_0\lt(f_x(s,y(s)\kappa_s(\T_s), v(s))y_3(s)+\frac12f_{xx}(s,y(s)\kappa_s(\T_s), v(s))y_1^2(s)\kappa_s(\T_s)\rt)\d s\rt]\\
&-\E\lt[\int^T_0\lt(f(s,y(s)\kappa_s(\T_s), v^\ep(s))-f(s,y(s)\kappa_s(\T_s), v(s))\rt)\kappa_s^{-1}(\T_s)\d s\rt]\\
&-\E\lt[\Phi_x(y(T)\kappa_T(\T_T))y_3(T)+\frac12\Phi_{xx}(y(T)\kappa_T(\T_T))y_1^2(T)\kappa_T(\T_T)\rt]
\le o(\ep).
\end{aligned}
\end{equation}
\end{lemma}
The proof of this lemma is given in the Appendix.
\medskip

For a pair of processes $(\varphi(\cdot),\psi(\cdot))$ in $L^2_\mathbb{H}(\Omega\times[0,T])\times L^2_\mathbb{H}(\Omega\times[0,T])$, we consider the following stochastic system:
\begin{equation}\label{eq_z}
\lt\{
\begin{array}{l}
\d z(t)=\Big(b_x(t,y(t)\kappa_t(\T_t),v(t))z(t)+\varphi(t)\Big)\d t +\Big(\beta_x(t,y(t)\kappa_t(\T_t),v(t))z(t)+\psi(t)\Big)\d W(t),\\
z(0)=0.
\end{array}\rt.
\end{equation}
With the help of this equation we define a linear functional
\begin{equation}\label{eq_phipsi}
I(\varphi(\cdot),\psi(\cdot))=\E\lt[\int^T_0f_x(t,y(t)\kappa_t(\T_t),v(t))z(t)\d t+\Phi_x(y(T)\kappa_T(\T_T))z(T)\rt],
\end{equation}
which is continuous in $L^2_\mathbb{H}(\Omega\times[0,T])\times L^2_\mathbb{H}(\Omega\times[0,T])$. The Riesz representation theorem yields that
there exists a unique pair of processes $(p(\cdot),K(\cdot))\in L^2_\mathbb{H}(\Omega\times[0,T])\times L^2_\mathbb{H}(\Omega\times[0,T])$ such that
\begin{equation}\label{eq_phipsi1}
I(\varphi(\cdot),\psi(\cdot))=\E\lt[\int^T_0(p(t)\varphi(t)+K(t)\psi(t))\d t\rt], \ (\varphi,\psi)\in L^2_\mathbb{H}(\Omega\times[0,T])\times L^2_\mathbb{H}(\Omega\times[0,T]).
\end{equation}
Notice that the processes of $p$ and $K$ do not depend on $(\varphi(\cdot),\psi(\cdot))$. By applying the above
representation result to the definition of $y_1$ in (\ref{eq_y_1}), we get that
\begin{equation}\label{eq_r1}
\begin{aligned}
&\E\lt[\int^T_0 f_x(s,y(s)\kappa_s(\T_s),v(s))y_1(s)\d s+ \Phi_x(y(T)\kappa_T(\T_T))y_1(T)\rt]\\
=&E\lt[\int^T_0 p(s)\underbrace{\Big(b(s,y(s)\kappa_s(\T_s),v^\ep(s))-b(s,y(s)\kappa_s(\T_s),v(s))\Big)\kappa_s^{-1}(\T_s)}_{=\varphi(s)}\d s\rt]\\
&+\E\lt[\int^T_0 K(s)\underbrace{\Big(\beta(s,y(s)\kappa_s(\T_s),v^\ep(s))-\beta(s,y(s)\kappa_s(\T_s),v(s))\Big)\kappa_s^{-1}(\T_s)}_{=\psi(s)}\d s\rt],
\end{aligned}
\end{equation}
and with the above indicated choice of $\varphi$ and $\psi$, $z$ defined by (\ref{eq_z}) coincides with $y_1$ defined by (\ref{eq_y_1}). Thus, (\ref{eq_phipsi}) and (\ref{eq_phipsi1}) yield (\ref{eq_r1}). With similar argument applied to (\ref{eq_y_2}), we have
\begin{equation}\label{eq_r2}
\begin{aligned}
&\E\lt[\int^T_0 f_x(s,y(s)\kappa_s(\T_s),v(s))y_2(s)\d s+ \Phi_x(y(T)\kappa_T(\T_T))y_2(T)\rt]\\
=&E\lt[\frac12\int^T_0\Big(p(s)b_{xx}(s,y(s)\kappa_s(\T_s),v(s))+K(s)\beta_{xx}(s,y(s)\kappa_s(\T_s),v(s))\Big) y_1^2(s)\kappa_s(\T_s)\d s\rt]\\
&+\E\lt[\int^T_0 p(s)\Big(b_x(s,y(s)\kappa_s(\T_s),v^\ep(s))-b_x(s,y(s)\kappa_s(\T_s),v(s))\Big)y_1(s)\d s\rt]\\
&+\E\lt[\int^T_0 K(s)\Big(\beta_x(s,y(s)\kappa_s(\T_s),v^\ep(s))-\beta_x(s,y(s)\kappa_s(\T_s),v(s))\Big)y_1(s)\d s\rt].
\end{aligned}
\end{equation}
We define a new (random) function $H$ by putting
$$
H(s,x,v,p,K)=\Big(f(s,x\kappa_s(\T_s),v)+pb(s,x\kappa_s(\T_s),v)+K\beta(s,x\kappa_s(\T_s),v)\Big)\kappa_s^{-1}(\T_s).
$$
Then, using (\ref{eq_r1}) and (\ref{eq_r2}), we can rewrite inequality (\ref{eq_var}) as
\begin{equation}\label{eq_var1}
\begin{aligned}
&-\E\lt[\int^T_0\Big(H(s,y(s),v^\ep(s),p(s),K(s))-H(s,y(s),v(s),p(s),K(s))\Big)\d s\rt]\\
&-\E\lt[\frac12\int^T_0H_{xx}(s,y(s),v(s),p(s),K(s))y^2_1(s)\d s +\frac12\Phi_{xx}\big(y(T)\kappa_T(\T_T)\big)y_1^2(T)\kappa_T(\T_T)\rt]\\
\le&o(\ep).
\end{aligned}
\end{equation}

Now we deal with the quadratic term. Let $Y(s):=y_1^2(s)$. Applying the It\^o formula to $Y(s)$, we get
\begin{equation}\label{eq:Y}
\begin{aligned}
\d Y(s)= &\lt[2Y(s)\Big(b_x(s,y(s)\kappa_s(\T_s),v(s))+\frac12\beta^2_x(s,y(s)\kappa_s(\T_s),v(s))\Big)+\Xi^\ep(s)\rt]\d s\\
&+\Big[2Y(s)\beta_x(s,y(s)\kappa_s(\T_s),v(s))+\Psi^\ep(s)\Big]\d W(s),
\end{aligned}
\end{equation}
where
\begin{equation*}
\begin{aligned}
\Xi^\ep(s)=&2 y_1(s)\ep^{-1}_s(\T_s)\Big(b(s,y(s)\kappa_s(\T_s),v^\ep(s))-b(s,y(s)\kappa_s(\T_s),v(s))\\
&+\beta_x(s,y(s)\kappa_s(\T_s),v(s))\big(\beta(s,y(s)\kappa_s(\T_s),v^\ep(s))-\beta(s,y(s)\kappa_s(\T_s),v(s))\big)\Big)\\
&+\ep^{-2}_s(\T_s)\lt(\beta(s,y(s)\kappa_s(\T_s),v^\ep(s))-\beta(s,y(s)\kappa_s(\T_s),v(s))\rt)^2
\end{aligned}
\end{equation*}
and
\begin{equation*}
\begin{aligned}
\Psi^\ep(s)=2 y_1(s)\kappa_s^{-1}(\T_s)\big(\beta(s,y(s)\kappa_s(\T_s),v^\ep(s))-\beta(s,y(s)\kappa_s(\T_s),v(s))\big).
\end{aligned}
\end{equation*}
For any $(\Xi(\cdot),\Psi(\cdot))$ in $L^2_\mathbb{H}(\Omega\times[0,T])\times L^2_\mathbb{H}(\Omega\times[0,T])$, we consider the following stochastic system:
\begin{equation}\label{eq:Z}
\lt\{
\begin{array}{l}
\d Z(s)=\lt[2Z(s)\Big(b_x(s,y(s)\kappa_s(\T_s),v(s))+\beta^2_x(s,y(s)\kappa_s(\T_s),v(s))\Big)+\Xi(s)\rt]\d s\\
\qquad\qquad+\Big[2Z(s)\beta_x(s,y(s)\kappa_s(\T_s),v(s))+\Psi(s)\Big]\d W(s),\\
Z(0)=0.
\end{array}
\rt.
\end{equation}
We define a new linear functional
\begin{equation}\label{eq:L}
L(\Xi(\cdot),\Psi(\cdot))=\E\lt[\int^T_0Z(s)H_{xx}(s,y(s),v(s),p(s),K(s))\d s+Z(T)\Phi_{xx}\big(y(T)\kappa_T(\T_T)\big)\kappa_T(\T_T)\rt],
\end{equation}
which too is continuous on $L^2_\mathbb{H}(\Omega\times[0,T])\times L^2_\mathbb{H}(\Omega\times[0,T])$. Using the same argument as above we see that there exists a unique pair of
$(P(\cdot), Q(\cdot))$ in $L^2_\mathbb{H}(\Omega\times[0,T])\times L^2_\mathbb{H}(\Omega\times[0,T])$ such that
\begin{equation}\label{eq:L1}
L(\Xi(\cdot),\Psi(\cdot))=\E\lt[\int^T_0\Big(P(s)\Xi(s)+Q(s)\Psi(s)\Big)\d s\rt].
\end{equation}
Now we apply the above result with $Z(s):=Y(s)=y_1^2(s), \Xi(s):=\Xi^\ep(s), \Psi(s):=\Psi^\ep(s)$, $s\in[0,T]$.
By using the estimates in Lemma \ref{lemma_y1y2}, we obtain that, for $1<p,q<\infty$ with $\frac1p+\frac1q=1$,
\[
\begin{aligned}
&\E\lt[\int^T_0y_1(s)\kappa_s^{-1}(\T_s)\big(b(s,y(s)\kappa_s(\T_s),v^\ep(s))-b(s,y(s)\kappa_s(\T_s),v(s))\big)P(s)\d s\rt]
\\
\le &C\E\lt[\sup_{s\in[0,T]}|y_1(s)|\sup_{s\in[0,T]}(\kappa_s^{-1}(\mathcal{T}_s))\ep^{1/2}\lt(\int^T_0|P(s)|^2 1_{[\tau-\ep,\tau+\ep]}(s)\d s\rt)^{1/2} \rt]\\
\le& C\lt(\E\lt[\sup_{s\in[0,T]}|y_1(s)|^{2p}\rt]\rt)^{1/2p} \lt(\E\lt[\sup_{s\in[0,T]}|\kappa_s(\mathcal{T}_s)|^{-2q}\rt]\rt)^{1/2q} \ep^{1/2}\lt(\E\lt[\lt(\int^T_0|P(s)|^2 1_{[\tau-\ep,\tau+\ep]}(s)\d s\rt)\rt]\rt)^{1/2}\\
\le& C\ep h(\ep),
\end{aligned}
\]
where, from the Dominated Convergence Theorem
\[
h(\ep)= \lt(\E\lt[\lt(\int^T_0|P(s)|^2 1_{[\tau-\ep,\tau+\ep]}(s)\d s\rt)\rt]\rt)^{1/2}\to 0, \ \ \textrm{when}\  \ep\to 0.
\]
Hence, we have
\[
\E\lt[\int^T_0y_1(s)\kappa_s^{-1}(\T_s)\big(b(s,y(s)\kappa_s(\T_s),v^\ep(s))-b(s,y(s)\kappa_s(\T_s),v(s))\big)P(s)\d s\rt]=o(\ep).
\]
Similarly, we get
\[
\E\lt[\int^T_02 y_1(s)\kappa_s^{-1}(\T_s)\big(\beta(s,y(s)\kappa_s(\T_s),v^\ep(s))-\beta(s,y(s)\kappa_s(\T_s),v(s))\big)Q(s)\d s\rt]=o(\ep).
\]
Therefore, the relations (\ref{eq:Y}), (\ref{eq:Z}), (\ref{eq:L}) and (\ref{eq:L1}) allow to rewrite inequality (\ref{eq_var1}) as
\begin{equation}\label{eq_var2}
\begin{aligned}
&-\E\lt[\int^T_0\Big(H(s,y(s),v^\ep(s),p(s),K(s))-H(s,y(s),v(s),p(s),K(s))\Big)\d s\rt]\\
&-\frac12\E\lt[\int^T_0 \kappa_s^{-2}(\T_s)\Big(\beta(s,y(s)\kappa_s(\T_s),v^\ep(s))-\beta(s,y(s)\kappa_s(\T_s),v(s))\Big)^2P(s)\d s\rt]\le o(\ep).
\end{aligned}
\end{equation}
Hence, by letting $\ep$ tend to zero, we deduce that
\begin{equation}\label{eq_varineq}
\begin{aligned}
&H(\tau,y(\tau),v,p(\tau),K(\tau))-H(\tau,y(\tau),v(\tau),p(\tau),K(\tau))\\
&+\frac12\kappa_\tau^{-2}(\T_\tau)\Big(\beta(\tau,y(\tau)\kappa_\tau(\T_\tau),v) -\beta(\tau,y(\tau)\kappa_\tau(\T_\tau),v(\tau))\Big)^2P(\tau)
\ge0
\end{aligned}
\end{equation}
holds for any $U$-valued $\F^B_\tau$-measurable random variable $v$, $\d \tau$-a.e., a.s.,
where we recall
$$
H(s,x,v,p,K)=\Big(f(s,x\kappa_s(\T_s),v)+pb(s,x\kappa_s(\T_s),v)+K\beta(s,x\kappa_s(\T_s),v)\Big)\kappa_s^{-1}(\T_s).
$$
Inequality (\ref{eq_varineq}) is the stochastic variational inequality of our control problem.
Since in our case the variational inequality is different from those in Peng \cite{P90} and Buckdahn \textit{et al.} \cite{BuDjLi}, we prefer to give a detailed proof of deriving (\ref{eq_varineq}) from (\ref{eq_var2}) in the Appendix.

\bigskip

Following similar arguments as the classical results of Bensoussan \cite{Be} and Peng \cite{P90}, the pair of processes $(p(\cdot),K(\cdot))$ is determined by an adjoint backward stochastic differential equation, i.e., $(p(\cdot),K(\cdot))$ is the unique solution of
\begin{equation}\label{eq_BSDE1}
\lt\{
\begin{array}{rl}
-\d p(s)&=\big[b_x(s,y(s)\kappa_s(\T_s),v(s))p(s)+\beta_x(s,y(s)\kappa_s(\T_s),v(s))K(s)+f_x(s,y(s)\kappa_s(\T_s),v(s))\big]\d s\\
&\qquad
-K(s)\d W(s)-\d N(s),\ \ s\in[0,T],\\
p(T)&=\Phi_x(y(T)\kappa_T(\T_T)),
\end{array}
\rt.
\end{equation}
and $(P(\cdot),Q(\cdot))$  is the unique solution of the following adjoint backward stochastic differential equation:
\begin{equation}\label{eq_BSDE2}
\lt\{
\begin{array}{rl}
-\d P(s)&=\big[2b_x(s,y(s)\kappa_s(\T_s),v(s))P(s)+\beta_x^2(s,y(s)\kappa_s(\T_s),v(s))P(s)\\
&\qquad\qquad+2\beta_x(s,y(s)\kappa_s(\T_s),v(s))Q(s)+H_{xx}(s,y(s),v(s),p(s),K(s))\big]\d s\\
&\qquad-Q(s)\d W(s)-\d M(s), \ \ s\in[0,T],\\
P(T)&=\Phi_{xx}(y(T)\kappa_T(\T_T))\kappa_T(\T_T),
\end{array}
\rt.
\end{equation}
where $N(\cdot)$ and $M(\cdot)$ are $\mathbb{H}$-adapted square integrable martingales orthogonal to $W$.
One can easily verify that the solutions $(p(\cdot),K(\cdot))$ and $(P(\cdot),Q(\cdot))$ satisfy  (\ref{eq_r1}) and (\ref{eq_r2}).

\begin{remark}
The two martingales $N(\cdot)$ and $M(\cdot)$ are introduced here from the Galtchouk-Kunita-Watanabe decomposition (we refer to \cite{KuWa}); this allows to guarantee the adaptedness of $(p(\cdot),K(\cdot))$ and $(P(\cdot),$ $Q(\cdot))$ with respect to $\mathbb{H}$. Such backward stochastic differential equations with respect to a non-Brownian filtration have been well studied, and they were also employed to study control problems, for instance, in Buckdahn and Ichihara \cite{BuIc} and Buckdahn \textit{et al.} \cite{BuLaRaTa}.
\end{remark}

As a consequence, we obtain the maximum principle theorem.
\begin{theorem}\label{thm_max}
Let (H1) and (H2) hold. If $(y(\cdot),v(\cdot))$ is the optimal solution of the control problem (\ref{eq_zeta}) and (\ref{eq_cost}), then we have
$$
(p(\cdot),K(\cdot))\in L^2_{\mathbb{H}}(\Omega\times[0,T])\times L^2_{\mathbb{H}}(\Omega\times[0,T])\ \ \textrm{and}\ \ (P(\cdot),Q(\cdot))\in L^2_{\mathbb{H}}(\Omega\times[0,T])\times L^2_{\mathbb{H}}(\Omega\times[0,T]),
$$
are the solutions of backward stochastic differential equations (\ref{eq_BSDE1}) and (\ref{eq_BSDE2}) respectively, such that the (stochastic) variational inequality (\ref{eq_varineq}) holds.
\end{theorem}

\section{Comparison with the Classical Case and Conclusion}
In this part we compare our result with the classical case, i.e., Peng's result \cite{P90}.

First, if $\sigma\equiv 0$, i.e.,if there is no fractional Brownian motion part, then obviously our result reduces to Peng's. Second, if $\sigma\neq 0$ but $H=1/2$, i.e., the fractional Brownian motion $B^H$ is nothing else but a  standard Brownian motion $B$, we show that our result coincides with Peng's characterisation of the optimal control. Here we only show that from equation (19) in Peng \cite{P90} we can obtain $(\ref{eq_BSDE1})$. With our notations, equation (19) in Peng \cite{P90} becomes
\begin{equation}\label{eq_p19}
\lt\{
\begin{array}{rl}
-\d p(s)&=\big[b_x(s,X^u(s),u(s))p(s)+\beta_x(s,X^u(s),u(s))K(s)+\sigma(s)K_1(s)+f_x(s,X^u(s),u(s))\big]\d s\\
&\qquad
-K(s)\d W(s)-K_1(s)\d B(s),\ \ s\in[0,T],\\
p(T)&=\Phi_x(X^u(T)).
\end{array}
\rt.
\end{equation}
We notice that in the classical case $H=1/2$, (\ref{ep}) yields $\kappa_s=\exp\lt\{\int^s_0\sigma(r)\d B(r)-\frac12\int^s_0\sigma^2(r)\d r\rt\}$. We put $\overline{p}(s):=p(s,\T_s)$, $\overline{K}_1(s)=K_1(s,\T_s)$ and $\overline{K}(s)=K(s,\T_s)$. Applying now standard arguments as above, and recalling the definition of $(y(\cdot,), v(\cdot))$ through (\ref{eq_y}),
we deduce that
\begin{equation}\label{eq_pbar}
\lt\{
\begin{array}{rl}
-\d \overline{p}(s)&=\big[b_x(s,y(s)\kappa_s(\T_s),v(s))\overline{p}(s)+\beta_x(s,y(s)\kappa_s(\T_s),v(s))\overline{K}(s)
+f_x(s,y(s)\kappa_s(\T_s),v(s))\big]\d s\\
&\qquad
-\overline{K}(s)\d W(s)-\overline{K}_1(s)\d B(s),\ \ s\in[0,T],\\
\overline{p}(T)&=\Phi_x(y(T)\kappa_T(\T_T)).
\end{array}
\rt.
\end{equation}
Equation (\ref{eq_pbar}) coincides with equation (\ref{eq_BSDE1}) with $N(t)=\int^t_0\overline{K}_1(s)\d B(s)$.
Hence our result is really a generalization of the classical one.

\section{Appendix}
In the Appendix we state some proofs of the results in Section 3.
\medskip

\noindent\textit{Proof of Theorem \ref{thm_girsanov}:} The proof is a bit technical and we split it into 3 steps.

\textbf{Step 1}: First we prove the existence and uniqueness of $\zeta^u$ in $L^p_{\mathbb{H}}(\Omega\times[0,T]),$ $p\ge 2$.  To this end, we define $\hat{\beta}$ and $\hat{b}$ as
$$
\hat{\beta}(t,x,u):=\kappa^{-1}_t(\mathcal{T}_t)\beta(t,x\kappa_t(\mathcal{T}_t),u),\ \
\hat{b}(t,x,u):=\kappa^{-1}_t(\mathcal{T}_t)b(t,x\kappa_t(\mathcal{T}_t),u),\ \ (t,x,u)\in[0,T]\times\R\times U.
$$
Then, from \textbf{(H2)} and (\ref{eq:ep}), $\hat{\beta}(\cdot,0,0)$ and $\hat{b}(\cdot,0,0) \in L^p_{\mathbb{H}}(\Omega\times[0,T]),$ $p\ge 2$. Furthermore, from \textbf{(H2)}, we know that there exists a constant $C>0$, such that
$$
|\hat{\beta}(t,x_1,u)-\hat{\beta}(t,x_2,u)|\le C|x_1-x_2|,\
|\hat{b}(t,x_1,u)-\hat{b}(t,x_2,u)|\le C|x_1-x_2|,\ x_1, x_2\in\R, (t,u)\in[0,T]\times U.
$$
Hence, equation (\ref{eq_zeta}) admits a unique solution $\zeta^u\in L^p_{\mathbb{H}}(\Omega\times[0,T])$, $p\ge 2$.
\medskip

\textbf{Step 2}:
Next we prove that $X^u$ is a solution of  equation (\ref{eq_sys}). Observe, that from the definition of $X^u$ and the above property of $\zeta^u$,  it follows that $X^u\in L^2_{\mathbb{H}}(\Omega\times[0,T])$. Let us
choose an arbitrary $F\in\mathcal{S}_{\K}$.  Then we have from (\ref{eq:gir}), for $t\in[0,T]$,
$$
\begin{aligned}
&\E\lt[F X^u(t)- F x_0\rt]=\E\lt[F(\T_t)\zeta^u(t)-F \zeta^u(0)\rt]\\
=&\E\lt[F(\T_t)x_0-Fx_0+F(\T_t)\int^t_0\kappa_s^{-1}(\mathcal{T}_s) \beta(s,\zeta^u(s)\kappa_s(\mathcal{T}_s),u(s,\mathcal{T}_s))\d W(s)\rt]\\ &+\E\lt[F(\T_t)\int^t_0\kappa_s^{-1}(\mathcal{T}_s)b(s,\zeta^u(s)\kappa_s(\mathcal{T}_s), u(s,\mathcal{T}_s))\d s\rt].
\end{aligned}
$$
We recall that, from (\ref{eq_defdeltaw}),
$$
\begin{aligned}
&\E\lt[F(\T_t)\int^t_0\kappa_s^{-1}(\mathcal{T}_s) \beta(s,\zeta^u(s)\kappa_s(\mathcal{T}_s),u(s,\mathcal{T}_s))\d W(s)\rt]\\
=&\E\lt[\int^t_0\kappa_s^{-1}(\mathcal{T}_s) \beta(s,\zeta^u(s)\kappa_s(\mathcal{T}_s),u(s,\mathcal{T}_s))D^W_s F(\T_t)\d s\rt].
\end{aligned}
$$
From the fact that $F\in \mathcal{S}_{\mathcal{K}}$ and the definition of $\mathcal{T}_t$, we deduce that (see, Jing and Le\'on \cite{JL} Page 7),
$$
\frac{\d F(\T_t)}{\d t}= \sigma(t)(\K^*\K D^BF)(\T_t,t).
$$
Using the above result, we obtain
\begin{equation}\label{eq_girsa}
\begin{aligned}
&\E\lt[F X^u(t)- F x_0\rt]\\
=&\E\lt[x_0\int^t_0\sigma(s)(\K^*\K D^BF)(\T_s,s)\d s\rt]\\
&+\E\lt[\int^t_0\kappa_s^{-1}(\mathcal{T}_s) \beta(s,\zeta^u(s)\kappa_s(\mathcal{T}_s),u(s,\mathcal{T}_s))D^W_s F(\T_s)\d s\rt]\\
&+\E\lt[\int^t_0\lt(\int^t_sD^W_s(\sigma(r)(\K^*\K D^BF)(\T_r,r))\d r \rt) \kappa_s^{-1}(\mathcal{T}_s)\beta(s,\zeta^u(s)\kappa_s(\mathcal{T}_s),u(s,\mathcal{T}_s))\d s\rt]\\
&+\E\lt[\int^t_0\kappa_s^{-1}(\mathcal{T}_s)b(s,\zeta^u(s)\kappa_s(\mathcal{T}_s), u(s,\mathcal{T}_s))F(\T_s)\d s\rt]\\
&+\E\lt[\int^t_0\lt(\int^t_s\sigma(r)(\K^*\K D^BF)(\T_r,r)\d r\rt)\kappa_s^{-1}(\mathcal{T}_s) b(s,\zeta^u(s)\kappa_s(\mathcal{T}_s),u(s,\mathcal{T}_s))\d s\rt].\\
\end{aligned}
\end{equation}
By applying the Fubini theorem, we get
$$
\begin{aligned}
I_1=&\E\lt[\int^t_0\lt(\int^t_sD^W_s\lt(\sigma(r)(\K^*\K D^BF)(\T_r,r)\rt)\d r\rt)\kappa_s^{-1}(\mathcal{T}_s) \beta(s,\zeta^u(s)\kappa_s(\mathcal{T}_s),u(s,\mathcal{T}_s))\d s\rt]\\
=&\int^t_0\E\lt[\int^r_0D^W_s(\sigma(r)(\K^*\K D^BF)(\T_r,r))\kappa_s^{-1}(\mathcal{T}_s) \beta(s,\zeta^u(s)\kappa_s(\mathcal{T}_s),u(s,\mathcal{T}_s))\d s\rt]\d r.
\end{aligned}
$$
Thus, taking into account that $\sigma(r)(\K^*\K D^BF)(\T_r,r)\in \mathcal{S}_{\mathcal{K}}$, we conclude from Remark \ref{remark} that
$$
\begin{aligned}
I_1=&\int^t_0\E\lt[\sigma(r)(\K^*\K D^BF)(\T_r,r)\int^r_0\kappa_s^{-1}(\mathcal{T}_s) \beta(s,\zeta^u(s)\kappa_s(\mathcal{T}_s),u(s,\mathcal{T}_s))\d W(s)\rt]\d r\\
=&\E\lt[\int^t_0\sigma(r)(\K^*\K D^BF)(\T_r,r)\int^r_0\kappa_s^{-1}(\mathcal{T}_s) \beta(s,\zeta^u(s)\kappa_s(\mathcal{T}_s),u(s,\mathcal{T}_s))\d W(s)\d r\rt].
\end{aligned}
$$
Consequently, using the Fubini Theorem now also for the latter double integral in (\ref{eq_girsa}), we get
\begin{equation*}
\begin{aligned}
\E\lt[F X^u(t)- F x_0\rt]=&\E\bigg[\int^t_0\sigma(r)(\K^*\K D^BF)(\T_r,r) \bigg\{x_0+\int^r_0\kappa_s^{-1}(\mathcal{T}_s) \beta(s,\zeta^u(s)\kappa_s(\mathcal{T}_s),u(s,\mathcal{T}_s))\d W(s)\\
&\qquad\qquad\qquad+\int^r_0\kappa_s^{-1}(\mathcal{T}_s)b(s,\zeta^u(s)\kappa_s(\mathcal{T}_s), u(s,\mathcal{T}_s))\d s\bigg\}\bigg]\\
&+\E\lt[\int^t_0\kappa_s^{-1}(\mathcal{T}_s) \beta(s,\zeta^u(s)\kappa_s(\mathcal{T}_s),u(s,\mathcal{T}_s))D^W_s F(\T_s)\d s\rt]\\
&+\E\lt[\int^t_0\kappa_s^{-1}(\mathcal{T}_s)b(s,\zeta^u(s)\kappa_s(\mathcal{T}_s), u(s,\mathcal{T}_s))F(\T_s)\d s\rt].\\
\end{aligned}
\end{equation*}
Hence, from (\ref{eq_zeta}) and by applying the Girsanov Theorem again, we get
$$
\begin{aligned}
&\E\lt[F X^u(t)- F x_0\rt]\\
=&\E\lt[\int^t_0(\K^*\K D^BF)(\T_s,s)\sigma(s)\zeta^u(s)\d s \rt]+ \E\lt[\int^t_0\kappa_s^{-1}(\mathcal{T}_s) \beta(s,\zeta^u(s)\kappa_s(\mathcal{T}_s),u(s,\mathcal{T}_s))D^W_s F(\T_s)\d s\rt]\\
&+\E\lt[\int^t_0\kappa_s^{-1}(\mathcal{T}_s)b(s,\zeta^u(s)\kappa_s(\mathcal{T}_s), u(s,\mathcal{T}_s))F(\T_s)\d s\rt]\\
=&\E\lt[\int^t_0(\K^*\K D^BF)(s)\sigma(s)X^u(s)\d s \rt]+ \E\lt[\int^t_0 \beta(s,X^u(s),u(s))D^W_s F\d s\rt]\\
&+\E\lt[F\int^t_0 b(s,X^u(s), u(s))\d s\rt].
\end{aligned}
$$
Since $\beta(\cdot,X^u,u)I_{[0,t]}$ is $\mathbb{H}$-adapted and square integrable, its Skorohod integral with respect to $W$ is well defined and coincides with the It\^o integral. Thus, from (\ref{eq_defdeltaw}), we have
\[
\E\lt[\int^t_0 \beta(s,X^u(s),u(s))D^W_s F\d s\rt] =\E\lt[F\int^t_0 \beta(s,X^u(s),u(s))\d W(s)\rt].
\]
Consequently,
\[
\E\lt[\int^t_0(\K^*\K D^B F)(s)\sigma(s)X^u(s)\d s\rt]=\E\lt[F G(t)\rt],
\]
where
\[
G(t)=X^u(t)-\lt(x_0+\int^t_0\beta(s,X^u(s),u(s))\d W(s)+\int^t_0b(s,X^u(s),u(s))\d s\rt).
\]
Observing that $\sigma X^u \in L^2_{\mathbb{H}}(\Omega\times[0,T])$ and $G(t)\in L^2(\Omega)$, we conclude from (\ref{eq_defdeltab}) that
$\sigma X^u I_{[0,t]}\in Dom\ \delta^B$ and
\[
\int^t_0 \sigma(s)X^u(s)\d B^H(s)=\delta(\sigma X^u I_{[0,t]})=G(t).
\]
This proves $X^u\in L^2_{\mathbb{H}}(\Omega\times[0,T])$ is a solution of (\ref{eq_sys}).
\medskip

\textbf{Step 3:} Now we prove the uniqueness. Suppose $X^u\in L^2_{\mathbb{H}}(\Omega\times[0,T])$ is a solution of (\ref{eq_sys}) such that $\sigma X^u I_{[0,t]}\in \ Dom\ \delta^B$, for every $t\in[0,T]$. Define $\eta$ as
$\eta(t)=X^u(t,\T_t)\kappa_t^{-1}(\T_t)$, $t\in[0,T]$. Then we have $\eta\in L^q_\mathbb{H}(\Omega\times[0,T])$, $1<q<2.$
For any $F\in\mathcal{S}_{\mathcal{K}}$, we have
\[
\begin{aligned}
&\E\lt[F\eta(t)-FX^u_0\rt]=\E\lt[F(A_t)X^u_t-FX_0^u\rt]\\
=&\E\lt[F(\mathcal{A}_t)\lt\{X^u_0+\int^t_0\sigma(s)X^u(s)\d B(s)+\int^t_0b(s,X^u(s),u(s))\d s +\int^t_0\beta(s,X^u(s),u(s))\d W(s)\rt\}\rt]\\
&-\E\lt[FX^u_0\rt].
\end{aligned}
\]
From the fact that
\[
\frac{\d F(\mathcal{A}_t)}{\d t}= -\sigma(t)(\K^*\K D^BF)(\mathcal{A}_t,t),
\]
and applying the same method as in \textbf{Step 2}, we deduce
\[
\begin{aligned}
&\E\lt[F\eta(t)-Fx_0\rt]\\
=&\E\bigg[-x_0\int^t_0\sigma(s)(\K^*\K D^B F(\mathcal{A}_s))(s)\d s+\int^t_0(\K^*\K D^B F(\mathcal{A}_s))(s)\sigma(s)X^u(s)\d s\\
&\qquad-\int^t_0\int^r_0\sigma(r)\K^*\K D^B(\K^*\K D^B F(\mathcal{A}_r)(r))(s)\sigma(s) X^u(s)\d s\d r\\
&\qquad+\int^t_0F(\mathcal{A}_s)b(s,X^u(s),u(s))\d s-\int^t_0\int^r_0\sigma(r)\K^*\K D^B F(\mathcal{A}_r)(r)b(s,X^u(s),u(s))\d s \d r\\
&\qquad+\int^t_0(D^W F(\mathcal{A}_s))(s)\beta(s,X^u(s),u(s))\d s\\
&\qquad-\int^t_0\int^r_0 \sigma(r)D^W(\K^*\K D^B F(\mathcal{A}_r)(r))\beta(s,X^u(s),u(s))\d s\d r\bigg].
\end{aligned}
\]
Since $X^u$ is a solution of (\ref{eq_sys}), we derive that
\[
\E\lt[F\eta(t)-Fx_0\rt]=\E\lt[\int^t_0F(\mathcal{A}_s)b(s,X^u(s),u(s))\d s+\int^t_0(D^W F(\mathcal{A}_s))(s)\beta(s,X^u(s),u(s))\d s\rt].
\]
We apply again the Girsanov Theorem and  (\ref{eq_defdeltaw}). Then
\[
\begin{aligned}
&\E\lt[F\eta(t)-Fx_0\rt]\\
=&\E\lt[F\int^t_0b(s,\eta(s)\kappa_s(\T_s),u(s,\T_s))\kappa_s^{-1}(\T_s)\d s+F\int^t_0\beta(s,\eta(s)\kappa_s(\T_s),u(s,\T_s))\kappa_s^{-1}(\T_s)\d W(s)\rt].
\end{aligned}
\]
From the arbitrariness of $F\in\mathcal{S}_{\K}$, we get
\[
\eta(t)=x_0+\int^t_0 b(s,\eta(s)\kappa_s(\T_s),u(s,\T_s))\kappa_s^{-1}(\T_s)\d s+\int^t_0\beta(s,\eta(s)\kappa_s(\T_s),u(s,\T_s))\kappa_s^{-1}(\T_s)\d W(s),
\]
$t\in[0,T]$.
But the $\mathbb{H}$-adapted continuous solution of this equation is unique and standard estimates show that it belongs to $S^2_{\mathbb{H}}$. Hence, $\eta\in L^2_{\mathbb{H}}(\Omega\times[0,T])$ is a solution of (\ref{eq_zeta}). Since equation (\ref{eq_zeta}) admits a unique solution, we have proved the uniqueness.
$\Box$
\bigskip

Let us present now the

\noindent\textit{Proof of Lemma \ref{lemma_y-y3}}: In this proof, for simplicity of notations, we make the conventions that
$V^\ep(s): =(s,y(s)\kappa_s(\mathcal{T}_s),v^{\ep}(s))$ and $V(s):=(s,y(s)\kappa_s(\mathcal{T}_s),v(s))$.
 Putting
\begin{equation*}
I=\int^t_0b\big(s,(y(s)+y_3(s))\kappa_s(\mathcal{T}_s),v^{\ep}(s)\big)\kappa_s^{-1}(\T_s)\d s +\int^t_0\beta\big(s,(y(s)+y_3(s))\kappa_s(\mathcal{T}_s),v^{\ep}(s)\big)\kappa_s^{-1}(\T_s)\d W(s),
\end{equation*}
we have, from the Taylor expansion, that
\begin{equation*}
\begin{aligned}
I= & \int^t_0\bigg[b\big(V^\ep(s)\big)\kappa_s^{-1}(\T_s)+ b_x\big(V^{\ep}(s)\big)y_3(s)\\
&\qquad+\lt(\int^1_0\int^1_0\lambda b_{xx}\big(s,y(s)\kappa_s(\mathcal{T}_s)+\lambda\mu y_3(s) \kappa_s(\mathcal{T}_s), v^{\ep}(s)\big)\d\lambda\d \mu \rt)y^2_3(s)\kappa_s(\T_s)\bigg]\d s\\
&+\int^t_0\bigg[\beta\big(V^{\ep}(s)\big)\kappa_s^{-1}(\T_s)+ \beta_x\big(V^{\ep}(s)\big)y_3(s)\\
&\qquad\quad+\lt(\int^1_0\int^1_0\lambda \beta_{xx}\big(s,y(s)\kappa_s(\mathcal{T}_s)+\lambda\mu y_3(s) \kappa_s(\mathcal{T}_s), v^{\ep}(s)\big)\d\lambda\d \mu \rt)y^2_3(s)\kappa_s(\T_s)\bigg]\d W(s),\\
\end{aligned}
\end{equation*}
which can be rewritten as
\begin{eqnarray*}
I&=&\int^t_0b\big(V(s)\big)\kappa_s^{-1}(\T_s)\d s+ \int^t_0\beta\big(V(s)\big)\kappa_s^{-1}(\T_s) \d W(s)\\
&&+\int^t_0b_x\big(V(s)\big)y_3(s)\d s+ \int^t_0\beta_x\big(V(s)\big)y_3(s)\d W(s)\\
&&+\int^t_0\lt(b\big(V^{\ep}(s)\big) -b\big(V(s)\big)\rt)\kappa_s^{-1}(\T_s)\d s
+\int^t_0\lt(\beta\big(V^{\ep}(s)\big) -\beta\big(V(s)\big)\rt)\kappa_s^{-1}(\T_s)\d W(s)\\
&&+\frac12\int^t_0b_{xx}\lt(V^{\ep}(s)\rt)y_3^2(s)\kappa_s(\T_s)\d s
+\frac12\int^t_0\beta_{xx}\big(V^{\ep}(s)\big)y_3^2(s)\kappa_s(\T_s)\d W(s)\\
&&+\int^t_0\lt(b_x\big(V^{\ep}(s)\big) -b_x\big(V(s)\big)\rt)y_3(s)\d s
+\int^t_0\lt(\beta_x\big(V^{\ep}(s)\big) -\beta_x\big(V(s)\big)\rt)y_3(s)\d W(s)\\
&&+\int^t_0\bigg(\int^1_0\int^1_0\lambda \Big(b_{xx}\big(s,y(s)\kappa_s(\mathcal{T}_s)+\lambda\mu y_3(s) \kappa_s(\mathcal{T}_s), v^{\ep}(s)\big)-b_{xx}\big(V(s)\big)\Big)\d\lambda\d \mu \bigg)y^2_3(s)\kappa_s(\T_s)\d s\\
&&+\int^t_0\bigg(\int^1_0\int^1_0\lambda \Big(\beta_{xx}\big(s,y(s)\kappa_s(\mathcal{T}_s)+\lambda\mu y_3(s) \kappa_s(\mathcal{T}_s), v^{\ep}(s)\big)\\
&&\qquad\qquad\qquad\quad-\beta_{xx}\big(V(s)\big)\Big)\d\lambda\d \mu\bigg) y^2_3(s)\kappa_s(\T_s)\d W(s).
\end{eqnarray*}
Consequently, according to the definitions of $y_1$ and $y_2$, we get
\begin{equation}\label{eq_I}
I = y(t) + y_3(t) + x_0 +\int^t_0 G^\ep(s)\d s+\int^t_0 \Lambda^\ep(s)\d W(s),
\end{equation}
where
\begin{equation}
\begin{aligned}
G^\ep(s)=&\frac12b_{xx}\big(V^{\ep}(s)\big)\lt(y_2^2(s)+2y_1(s)y_2(s)\rt)\kappa_s(\T_s)
+\Big(b_x\big(V^{\ep}(s)\big) -b_x\big(V(s)\big)\Big)y_2(s)\\
&+\bigg(\int^1_0\int^1_0\lambda \Big(b_{xx}\big(s,y(s)\kappa_s(\mathcal{T}_s)+\lambda\mu y_3(s) \kappa_s(\mathcal{T}_s), v^{\ep}(s)\big)-b_{xx}\big(V(s)\big)\Big)\d\lambda\d \mu \bigg)y^2_3(s)\kappa_s(\T_s)
\end{aligned}
\end{equation}
and
\begin{equation}
\begin{aligned}
\Lambda^\ep(s)=&\frac12\beta_{xx}\big(V^{\ep}(s)\big)
\lt(y_2^2(s)+2y_1(s)y_2(s)\rt)\kappa_s(\T_s)+\Big(\beta_x\big(V^{\ep}(s)\big) -\beta_x\big(V(s)\big)\Big)y_2(s)\\
&+\bigg(\int^1_0\int^1_0\lambda \Big(\beta_{xx}\big(s,y(s)\kappa_s(\mathcal{T}_s)+\lambda\mu y_3(s) \kappa_s(\mathcal{T}_s), v^{\ep}(s)\big)-\beta_{xx}\big(V(s)\big)\Big)\d\lambda\d \mu\bigg) y^2_3(s)\kappa_s(\T_s).
\end{aligned}
\end{equation}
We consider now the estimate of $\sup_{0\le t\le T}\E\lt[\lt|\int^t_0\Lambda^\ep(s)\d W(s)\rt|^2\rt]$. From the Burkh\"older-Davis-Gundy inequality we have
\[
\begin{aligned}
&\sup_{0\le t\le T}\E\lt[\lt|\int^t_0\Lambda^\ep(s)\d W(s)\rt|^2\rt]\le C \E\lt[\int^T_0\lt|\Lambda^\ep(s)\rt|^2\d s\rt]\\
\le&C \E\bigg[\int^T_0\lt(y^2_2(s)+2y_1(s)y_2(s)\rt)^2\kappa_s^2(\T_s)\d s+\int^T_0\Big(\beta_x\big(V^{\ep}(s)\big) -\beta_x\big(V(s)\big)\Big)^2y_2^2(s)\d s\bigg]\\
&+\int^T_0 y^4_3(s)\kappa_s^2(\T_s)\lt(\int^1_0 \lt(\lt|\beta_{xx}\big(s,y(s)\rt|\kappa_s(\mathcal{T}_s)+\theta \lt|y_3(s)\rt| \kappa_s(\mathcal{T}_s), v^{\ep}(s)\big)+\lt|\beta_{xx}\big(V(s)\big)\rt|\rt)\d \theta\rt)^2\d s.
\end{aligned}
\]
By using the estimates in Lemma \ref{lemma_y1y2} and applying the Dominated Convergence Theorem, we obtain
$$
\sup_{0\le t\le T}\E\lt[\lt|\int^t_0\Lambda^\ep(s)\d W(s)\rt|^2\rt]=o(\ep^2).
$$
Similar arguments can be applied to estimate $\sup_{0\le t\le T}\E\lt[\lt|\int^t_0G^\ep(s)\d s\rt|^2\rt]$. Hence we have
$$
\sup_{0\le t\le T}\E\lt[\lt|\int^t_0G^\ep(s)\d s\rt|^2+\lt|\int^t_0\Lambda^\ep(s)\d W(s)\rt|^2\rt]=o(\ep^2).
$$
Therefore, from (\ref{eq_I}) and the definition of $y^\ep$, we get that
$$
\begin{aligned}
y^\ep(t)-y(t)-y_3(t)
=&\int^t_0A^\ep(s)(y^\ep(s)-y(s)-y_3(s))\d s+\int^t_0\Theta^\ep(s)(y^\ep(s)-y(s)-y_3(s))\d W(s)\\
&+\int^t_0G^\ep(s)\d s+\int^t_0\Lambda^\ep(s)\d W(s),
\end{aligned}
$$
with the both factors
$$
A^\ep(s)=b_x(s,(y^\ep(s)-\theta(y(s)+y_3(s))),v^\ep(s)),\ \theta\in[0,1],
$$
and
$$
\Theta^\ep(s)=\beta_x(s,(y^\ep(s)-\lambda(y(s)+y_3(s))),v^\ep(s)),\ \lambda\in[0,1],
$$
which are being uniformly bounded according to \textbf{(H2)}. Finally, we can derive our estimate by applying standard arguments. $\Box$
\bigskip

\noindent\textit{Proof of Lemma \ref{lemma_var}:} From the optimality of $(y(\cdot),v(\cdot))$, we have
\begin{equation}
\begin{aligned}
0\le J(v^\ep)-J(v)=&\E\lt[\int^T_0f(s,y^\ep(s)\kappa_s(\T_s),v^\ep(s))\ep^{-1}_s(\T_s)\d s +\Phi(y^\ep(T)\kappa_T(\T_T))\kappa_T^{-1}(\T_T)\rt]\\
-&\E\lt[\int^T_0f(s,y(s)\kappa_s(\T_s),v(s))\ep^{-1}_s(\T_s)\d s +\Phi(y(T)\kappa_T(\T_T))\kappa_T^{-1}(\T_T)\rt].
\end{aligned}
\end{equation}
The Lemmata \ref{lemma_y1y2} and \ref{lemma_y-y3} lead to
\begin{eqnarray*}
0&\le&\E\lt[\int^T_0\Big(f(s,(y(s)+y_3(s))\kappa_s(\T_s),v^\ep(s))-f(s,y(s)\kappa_s(\T_s),v(s))\Big) \kappa_s^{-1}(\T_s)\d s\rt]\\
&&+\E\lt[\Big(\Phi((y(T)+y_3(T))\kappa_T(\T_T))-\Phi(y(T)\kappa_T(\T_T))\Big)\kappa_T^{-1}(\T_T)\rt]+o(\ep)\\
&=&\E\lt[\int^T_0\Big(f(s,(y(s)+y_3(s))\kappa_s(\T_s),v(s))-f(s,y(s)\kappa_s(\T_s),v(s))\Big) \kappa_s^{-1}(\T_s)\d s\rt]\\
&&+\E\lt[\Big(\Phi((y(T)+y_3(T))\kappa_T(\T_T))-\Phi(y(T)\kappa_T(\T_T))\Big)\kappa_T^{-1}(\T_T)\d s\rt]\\
&&+\E\lt[\int^T_0\Big(f(s,(y(s)+y_3(s))\kappa_s(\T_s),v^\ep(s))-f(s,(y(s)+y_3(s))\kappa_s(\T_s),v(s))\Big) \kappa_s^{-1}(\T_s)\d s\rt]+o(\ep).
\end{eqnarray*}
Hence, by applying Taylor's expansion up to the second order, we obtain
\begin{eqnarray*}
0&\le&\E\lt[\int^T_0\Big(f_x(s,y(s)\kappa_s(\T_s),v(s))y_3(s)+\frac12f_{xx}(s,y(s)\kappa_s(\T_s),v(s))\Big)y_3^2(s) \kappa_s(\T_s)\d s\rt]\\
&&+\E\lt[\int^T_0\Big(f(s,y(s)\kappa_s(\T_s),v^\ep(s))-f(s,y(s)\kappa_s(\T_s),v(s))\Big) \kappa_s^{-1}(\T_s)\d s\rt]\\
&&+\E\lt[\int^T_0\Big(f_x(s,y(s)\kappa_s(\T_s),v^\ep(s))-f_x(s,y(s)\kappa_s(\T_s),v(s))\Big) y_3(s)\d s\rt]\\
&&+\E\lt[\frac12\int^T_0\Big(f_{xx}(s,y(s)\kappa_s(\T_s),v^\ep(s))-f_{xx}(s,y(s)\kappa_s(\T_s),v(s))\Big) y_3^2(s)\kappa_s(\T_s)\d s\rt]\\
&&+\E\lt[\Phi_x(y(T)\kappa_T(\T_T))y_3(T)+\frac12\Phi_{xx}(y(T)\kappa_T(\T_T))y_3^2(T)\kappa_T(\T_T)\rt]+o(\ep).
\end{eqnarray*}
The desired inequality is obtained by using the hypothesis (H2) and Lemma \ref{lemma_y1y2}. $\Box$
\medskip

\textit{Proof $($of $($\ref{eq_varineq}$)$$)$:} For any $\tilde{v}\in \mathcal{U}_{ad}$, we define a new admissible control
$$
\tilde{v}^\ep(t)=
\lt\{
\begin{array}{ll}
\tilde{v}(t),& t\in[\tau-\ep,\tau+\ep];\\
v(t),& t\in[0,T]\backslash [\tau-\ep,\tau+\ep].
\end{array}
\rt.
$$
Let us put
$$\begin{aligned}
\Theta^{\tilde{v}}(s):=&H(s,y(s),\tilde{v}^\ep(s),p(s),K(s))-H(s,y(s),v(s),p(s),K(s))\\
&+\frac12\kappa^{-2}_s(\T_s)\Big(\beta(s,y(s)\kappa_s(\T_s),\tilde{v}^{\ep}(s)) -\beta(s,y(s)\kappa_s(\T_s),v(s))\Big)^2P(s),
\end{aligned}
$$
and assume that (\ref{eq_varineq}) does not hold. Then there exist $\delta>0$ and an admissible control $\tilde{v}$ such that
the set $\Lambda^{\tilde{v}}:=\{(s,\omega):\Theta^{\tilde{v}}(s)(\omega)\le -\delta\}$ satisfies
\begin{equation}\label{eq_v1}
\E\lt[\lambda\lt(\Lambda^{\tilde{v}}\rt)\rt]\ge \delta>0,
\end{equation}
where $\lambda\lt(\Lambda^{\tilde{v}}\rt)=\int^T_0 I_{\Lambda^{\tilde{v}}}(s)\d s$. We derive from (\ref{eq_v1}) that
\begin{equation}\label{eq_v2}
\E\lt[\lambda\lt(\Lambda^{\tilde{v}}\cap [0,T/2]\rt)\rt]\ge \delta/2\ \ \textrm{or}\ \
\E\lt[\lambda\lt(\Lambda^{\tilde{v}}\cap [T/2,T]\rt)\rt]\ge \delta/2.
\end{equation}
Hence there exists $\tau_1\in[0,T/2]$ such that
\begin{equation}\label{eq_v3}
\E\lt[\lambda\lt(\Lambda^{\tilde{v}}\cap [\tau_1,\tau_1+ T/2]\rt)\rt]\ge \delta/2.
\end{equation}
Similarly, from (\ref{eq_v3}), we get that there exists $\tau_2\in[\tau_1,\tau_1+T/4]$ such that
\begin{equation}\label{eq_v4}
\E\lt[\lambda\lt(\Lambda^{\tilde{v}}\cap [\tau_2,\tau_2+ T/4]\rt)\rt]\ge \delta/4,
\end{equation}
etc. Consequently, for $n\ge 2$, there exists $\tau_n\in[\tau_{n-1},\tau_{n-1}+T/2^n]$ such that
\begin{equation}\label{eq_vn}
\E\lt[\lambda\lt(\Lambda^{\tilde{v}}\cap [\tau_n,\tau_n+ T/2^n]\rt)\rt]\ge \delta/2^{n}.
\end{equation}
Furthermore, there exists $\tau\in[0,T]$ with $\tau_n\to \tau$ $(n\to \infty)$, and $|\tau_n-\tau|\le T/2^n$, $n\ge 1.$ Hence, we have
\begin{equation}\label{eq_vnn}
\E\lt[\lambda\lt(\Lambda^{\tilde{v}}\cap [\tau-\ep_n,\tau+ \ep_n]\rt)\rt]\ge \E\lt[\lambda\lt(\Lambda^{\tilde{v}}\cap [\tau_n,\tau_n+ \ep_n/2]\rt)\rt] \ge \delta\ep_n/2T,
\end{equation}
where $\ep_n=2T/2^n$, $n\ge 1$.

We define
$$
v^{\ep_n}(t)=
\lt\{
\begin{array}{ll}
\overline{v}(t),& t\in[\tau-\ep_n,\tau+\ep_n];\\
v(t),& t\in[0,T]\backslash [\tau-\ep_n,\tau+\ep_n],
\end{array}
\rt.
$$
where $\overline{v}(t):=\tilde{v}(t)I_{\lt\{\Theta^{\tilde{v}}\le -\delta,\ t\in[\tau-\ep_n,\tau+\ep_n]\rt\} }+ v(t)I_{\lt\{\Theta^{\tilde{v}}> -\delta \ \textrm{or}\  t\notin [\tau-\ep_n,\tau+\ep_n]\rt\}}$. It follows that
\[
\Theta^{\overline{v}}(t)=\Theta^{\tilde{v}}(t)I_{\lt\{\Theta^{\tilde{v}}\le -\delta,\ t\in[\tau-\ep_n,\tau+\ep_n]\rt\} }(t).
\]
From (\ref{eq_var2}), we derive that
\[\begin{aligned}
o(\ep_n)\le& \E\lt[\int^{\tau+\ep_n}_{\tau-\ep_n}\Theta^{\overline{v}}(t)\d t\rt]
=\E\lt[\int^{\tau+\ep_n}_{\tau-\ep_n}\Theta^{\overline{v}}(t)I_{\lt\{\Theta^{\tilde{v}}\le -\delta\rt\}}\d t\rt]\\
\le & -\delta \E\lt[\lambda\lt(\Lambda^{\tilde{v}}\cap[\tau-\ep_n,\tau+\ep_n]\rt)\rt]\le -\delta \ep_n/2T, \ n\ge 1.
\end{aligned}
\]
This leads to contradiction. Consequently, $\Theta^{\tilde{v}}\ge 0$, a.s., $\d s$-a.e., for any $\tilde{v}\in\mathcal{U}_{ad}$, in particular, for $\tilde{v}\equiv v$, an $\mathcal{F}_{\tau}^B$-measurable random variable. $\Box$


\medskip


\begin{thebibliography}{99}
\bibitem{Be}Bensoussan A. Lecture on Stochastic Control, in Nonlinear Filtering and Stochastic Control, Lecture Notes in Mathematics 972, Springer-Verlag, 1981.
\bibitem{BHOS}Biagini F, Hu Y, {\O}ksendal B, Sulem A. A stochastic maximum principle for processes driven by fractional Brownian motion. Stochastic Processes and their Applications, 100 (2002) 233-253.
\bibitem{Bi} Bismut J M. An introductory approach to duality in optimal stochastic control. SIAM Review, 20 (1978), 62-78.
\bibitem{Bu} Buckdahn R. Anticipative Girsanov Transformations and Skorohod Stochastic Differential Equations, Memoirs of the AMS, 111, N.533, 1994.
\bibitem{BuDjLi} Buckdahn R, Djehiche B, Li J.  A general stochastic maximum principle for SDEs of mean-field type. Applied Mathematics $\&$ Optimization,  64 (2011), 197-216.
\bibitem{BuIc} Buckdahn R, Ichihara N. Limit theorem for controlled backward SDEs and homogenization of Hamilton-Jacobi-Bellman equations. Applied Mathematics $\&$ Optimization, 51 (2005), 1-33.
\bibitem{BuLaRaTa} Buckdahn R, Labed B, Rainer C, Tamer L. Existence of an optimal control for stochastic control systems with nonlinear cost functional. Stochastics: An International Journal of Probability and Stochastics Processes, 82 (2010), 241-256.
\bibitem{ChNu} Cheridito P, Nualart D. Stochastic integral of divergence type with respect to fractional Brownian motion with Hurst parameter $H\in(0,\frac12)$. Annales de l'Institut Henri Poincar\'e, 41 (2005), 1049-1081.
\bibitem{GuNu} Guerra J, Nualart D. Stochastic differential equations driven by fractional Brownian motion and standard Brownian motion. Stochastic Analysis and Applications,  26 (2008), 1053-1075.
\bibitem{HaHuSo} Han Y, Hu Y, Song J. Maximum principle for general controlled systems driven by fractional Brownian motions. Applied Mathematics $\&$ Optimization, 7 (2013), 279-322.
\bibitem{HuZh} Hu Y, Zhou X. Stochastic control for linear systems driven by fractional noises. SIAM J. Control Optim., 43 (2005), 2245-2277.
\bibitem{JiMa} Jien Y, Ma J. Stochastic differential equations driven by fractional Brownian motions. Bernoulli, 15 (2009), 846-870.
\bibitem{JL} Jing S, Le\'on J A. Semilinear backward doubly stochastic differential equations and SPDEs driven by fractional Brownian motion with Hurst parameter in (0,1/2). Bulletin des Sciences Math\'ematiques, 135 (2011), 896-935.
\bibitem{KuWa} Kunita H, Watanabe S. On square integrable martingales. Nagoya Mathematical Journal, 30 (1967), 209-245.
\bibitem{LeNu} Le\'on J A, Nualart D. An extension of the divergence operator for Gaussian processes. Stochastic Processes and Their Applications, 115 (2005), 481-492.
\bibitem{LeSa} Le\'on J A, San~Mart\'in J. Linear stochastic differential equations driven by a fractional Brownian motion with Hurst parameter less than 1/2. Stochastic Analysis and Applications, 25 (2007), 105-126.
\bibitem{MiSh} Mishura Y S, Shevchenko G M. Existence and uniqueness of the solution of stochastic differential equation involving Wiener process and fractional Brownian motion with Hurst index $H>1/2$. Communications in Statistics - Theory and Methods, 40 (2011), 3492-3508.
\bibitem{Nu} Nualart D. Stochastic integration with respect to fractional Brownian motion and applications. In: Stochastic Models, Proceedings of the Seventh Symposium on Probability and Stochastic Processes, ed. by J. M. Gonz\'alez-Barrios et al. Contemporary Mathematics, 336 (2003), 3-39.
\bibitem{P90} Peng S. A general stochastic maximum principle for optimal control problems. SIAM Jornal of Control and Optimization, 28, (1990), 966-979.
\end{thebibliography}
\end{document}